\documentclass[11pt]{article}
\usepackage{mathrsfs}
\usepackage{amssymb,latexsym}
\usepackage{amsmath}
\usepackage{graphicx}
\usepackage{setspace}
\usepackage{amsmath,amsfonts,amssymb,amsbsy,amsthm,latexsym,lscape}
\usepackage{color}

\def\oneandahalf{\par\baselineskip=18pt}

\def\E{{\rm E}}
	{\end{list}}
\newenvironment{hangref}{\begin{list}{}{\setlength{\itemsep}{0pt}%
			\setlength{\parsep}{0pt}\setlength{\leftmargin}{+\parindent}%
			\setlength{\itemindent}{-\parindent}}}{\end{list}}

\def\standard{\oddsidemargin=0in
	\evensidemargin=0in
	\topmargin =-.5in
	\textheight=9.0in
	\textwidth=6.5in}

\standard
\newtheorem{theorem}{Theorem}
\newtheorem{proposition}{Proposition}

\newtheorem{assumption}{Assumption}

\makeatletter
\def\thanks#1{\protected@xdef\@thanks{\@thanks
		\protect\footnotetext{#1}}}
\makeatother

\begin{document}

\normalsize

\title{
Some Studies on Stochastic Optimization based Quantitative Risk Management\thanks{The paper is published: Zhaolin Hu. 2025. Some studies on stochastic optimization based quantitative risk management. {\it Operations Research Transactions}, 29(3) 135-159.}}

\author{
        Zhaolin Hu\\
        {\footnotesize
        School of Economics and
        Management, Tongji University, Shanghai, 200092, China, russell@tongji.edu.cn}
}

\date{}

\maketitle

\begin{abstract}
Risk management often plays an important role in decision making under uncertainty. In quantitative risk management, assessing and optimizing risk metrics requires efficient computing techniques and reliable theoretical guarantees. In this paper, we introduce several topics on quantitative risk management and review some of the recent studies and advancements on the topics.  We consider several risk metrics and study decision models that involve the metrics, with a main focus on the related computing techniques and theoretical properties. We show that stochastic optimization, as a powerful tool, can be leveraged to effectively address these problems.
\end{abstract}

\oneandahalf

\section{Introduction}\label{sec:intro}

In decision making, there often exist profound risk and uncertainty. Risk management is an important component for the management practice. It is particularly important for the decision makers who are risk averse. Decision makers are concerned about the potential risk of an outcome and would like to assess how large is the potential loss. In management science, finance, engineering, and many other fields, the risk measures that quantify the potential risk of an outcome are introduced and widely studied. These quantities are often fundamental for institutions and individuals. For instance, an investment agency may need to set its capital reservation based on the risk exposed according to some regulation rule. An individual investor would like to know the risk of the chosen investment portfolio. When doing decision, decision makers often incorporate the risk metrics into some optimization models, to optimally balance the tradeoff between return and risk. To address these models, optimization techniques are developed and widely implemented. 

In this paper, we review some of the studies on stochastic optimization based risk management. We also conduct some new studies on these topics. As a unified approach, we treat many risk measures from the perspective of stochastic optimization. By doing so, we emphasize that, we can make use of the theories and methods of stochastic optimization. Especially, we can often deal with the nonlinearity of the risk measures and transform the problems of risk measures to the problems involving expectation. The expectation functional is linear and is often easier to handle. This bridge allows us to use tools for expectation to address a number of risk measures.   

Given a risk measure, how to assess its value is essential for management practice. This problem is referred to as estimation of risk. Due to the complexity of loss position and risk measure, evaluating the risk measure value of a loss position is often difficult. In risk management, simulation techniques are widely used. In the simulation techniques, sample points (scenarios) of the random variates are generated and the outcomes are evaluated for these sample points. The outcomes are then aggregated to conduct the evaluation or estimation of the metric of interest. Glasserman (2004) provided a systematic studies of using Monte Carlo methods to address various financial estimation problems. In risk management, the risk values may depend on various parameters. It is often beneficial to assess the change of risk with respect to the change of parameters. Such analysis is categorized as sensitivity analysis. Sensitivity analysis is an important problem in risk management and has been studied extensively in the past decades. Besides estimation and sensitivity analysis, the optimization of the risk measure is another essential problem and is widely studied in the literature, typically studied within the framework of stochastic optimization.  Ruszczy\'nski and Shapiro (2006) studied optimization of convex risk measures. They built various theoretical properties of the risk optimization problems. Over the past decades, stochastic optimization has become a powerful tool, especially in the era of data. Both theory and application of stochastic optimization have been studied extensively, see, e.g., Shapiro et al. (2014). In this paper, we use stochastic optimization theories and methods to address the computing problems of the risk measures.       

Among the various risk measures, the value-at-risk (VaR) and conditional value-at-risk (CVaR) received the most attentions. VaR is defined as the critical value that the random loss will not exceed with a given confidence level. Mathematically, it is the quantile of a random position, see, e.g., Jorion (2006) for a comprehensive introduction and study of VaR. CVaR, as an alternative of VaR, is defined as the average loss of a random position beyond the VaR, see, e.g., Rockafellar and Uryasev (2000). Hong et al. (2014) reviewed the Monte Carlo methods of VaR and CVaR. They provided a comprehensive introduction for the estimation of VaR and CVaR, the sensitivity analysis of VaR and CVaR, and the optimization of VaR and CVaR. While the two risk measures are widely used, they also expose some inefficiency in some practical applications. This raised a question about what kind of risk measure should be used to address specific risk management problems. The seminal paper Artzner et al. (1999) started the axiomatic study for risk measure. They introduced the notion of coherent risk measure. Shortly after, F\"{o}llmer and Schied (2002) and Frittelli and Gianin (2002) relaxed the condition of the coherent risk measure and instead proposed the notion of convex risk measure. In this paper, our focus is on some classes of important risk measures, especially the convex risk measure. We discuss the computation of several classes of risk measures. We notice that many risk measures can be expressed as the optimal value of a stochastic optimization problem involving expectation. We first consider a class of risk measures that can be defined via the optimization of an expectation function. We focus on the optimized certainty equivalent (OCE) risk measure. The OCE risk measure includes the well known CVaR as special instance. We then consider a class of risk measures that can be defined as a stochastic optimization problem with expectation constraint. We mainly analyze the utility-based shortfall risk (SR). The OCE risk measure and the SR are two classes of most widely used convex risk measures.         

We consider the estimation, the sensitivity analysis and the optimization of these risk measures. For the estimation of OCE risk measure, we implement the sample average approximation (SAA) to construct the estimator and build the consistency and asymptotic normality of the estimator. We derive the derivative sensitivity expression for the OCE risk measure, propose the corresponding estimator, and analyze the properties of the estimator. For the optimization of the OCE risk measure, we also apply the SAA method and build the statistical properties of the method. For the SR, we also consider the three computing problems. We review the main studies that addressed these problems.    

In this paper, we also consider an important risk measure, the probability. We focus on the optimization of this risk measure. The optimization problem with a probability constraint is well known as chance constrained program (CCP). From a mathematical programming perspective, when there exists randomness in the constraints of an optimization problem, it is natural to require the constraints be satisfied with a given high probability. This results in the formulation of CCP. The CCP was introduced in the 1950s by Charnes et al. (1958). It has become an important stochastic optimization model since then. In this paper, we review some of the studies for CCPs. 

The rest of this paper is organized as follows. In Section \ref{sec:oce}, we analyze the OCE risk measure and build some statistical results for the computing. In Section \ref{sec:SR}, we review the computing studies for SR. In Section \ref{sec:CCP}, we review various optimization approaches for CCP. We conclude the paper and make some further discussions in Section \ref{sec:discussion}.

\section{Optimized Certainty Equivalent Risk Measure}\label{sec:oce}

The notion of optimized certainty equivalent (OCE) was introduced by  Ben-Tal and Teboulle (1986, 2007). We
follow the notions in Ben-Tal and Teboulle (2007) and Hamm et al.
(2013) to formally define the OCE risk measure. Consider a financial
position $X$ with the distribution $P$. Let $u$ be a proper closed
concave nondecreasing utility function, which can be normalized to
satisfy $u(0)=0$ and $1\in \partial u(0)$ where $\partial$ denotes
the subdifferential. We first define the OCE of $X$ as follows
\begin{equation*}
\mbox{OCE}_u(X)=\mathop{\mbox{sup}}_{\eta\in \mathbb{R}}\ \left\{
\eta+\E\left[u(X-\eta)\right] \right\}.
\end{equation*}
OCE has a clear interpretation: given a position $X$, suppose the
decision maker consumes $\eta$ of $X$. Then the resulting present
value of $X$ is $\eta+\E\left[u(X-\eta)\right]$. Thus OCE is
actually the present value generated by an optimal allocation
between the present and future consumption (Ben-Tal and Teboulle
2007).

The OCE risk measure is then defined as the negative OCE
\begin{equation*}
\mbox{OCE-R}_u(X):=-\mbox{OCE}_u(X).
\end{equation*}
The OCE risk measure includes some famous risk measures. It is shown
in Ben-Tal and Teboulle (2007) that when the utility function $u(z)
= 1-e^{-z}$ the OCE risk measure is the entropic risk measure, i.e.
$\mbox{OCE-R}_u(X)= \log\E[e^{- X}]$. Ben-Tal and Teboulle (2007)
also built that an OCE risk measure is a coherent risk measure if
and only if the utility function $u(z)$ takes a piecewise linear
utility form as follows
\begin{equation*}\label{OCElinear}
u(z)=\gamma_1 \left[z\right]^+-\gamma_2 \left[-z\right]^+
\end{equation*}
with  $0\le \gamma_1<1<\gamma_2$. Setting $\gamma_1=0$ and
$\gamma_2=1/\alpha$, the utility function above yields the well
known condition value-at-risk (CVaR), i.e.,
\begin{equation*}
  \hbox{CVaR}_{1-\alpha}(X) =
  \inf_{t\in \mathbb{R}}\left\{t +
  \frac{1}{\alpha}\E[-X-t]^+\right\}.
\end{equation*}

\subsection{Estimation of OCE Risk Measure}\label{sec:estioce}

Following the convention, we consider the random loss $\xi = -X$.
Define $l(z)=-u(-z)$. Then we can naturally take the OCE risk
measure as the optimal value of the following stochastic
optimization problem
\begin{eqnarray}\label{estioce}
&\displaystyle\mathop{\mbox{minimize}}_{t\in T} & t +\E[l(\xi-t)],
\end{eqnarray}
where the feasible set $T$ is a compact interval that is specified reasonably (it
can be extended to $\mathbb{R}$ when necessary). In this paper we assume that such $T$ can always be specified. Throughout this paper we consistently use $t$ to denote the optimization variable that is used to generate the risk measures. Note that the 
definition of OCE risk measure represented by (\ref{estioce}) is directly proposed in some other 
studies, see, for instance, Rockafellar and
Uryasev (2013) and Rockafellar and Royset (2015). In these papers,
$l$ is called a regret.

Based on the expression above, Hamm et al. (2013) proposed a
stochastic root finding approach for estimating the OCE risk
measure. They used a two-step approach where the first step is to estimate the optimal solution $t^*$ of Problem (\ref{estioce}), and the second step is to plug the estimator of $t^*$ into the objective function in Problem (\ref{estioce}) and then to estimate the objective function. To estimate $t^*$, Hamm et al. (2013) applied the first order optimality condition for Problem (\ref{estioce}) and showed that with some technical conditions $t^*$ is the unique root of the equation 
\begin{equation}\label{equ:OCEroot}
1-\E[l'(\xi-t)]=0.        
\end{equation}
They used the Robbins-Monro algorithm to search the root of Equation(\ref{equ:OCEroot}) and built the statistical properties based on the study of Dunkel and Weber (2010), who developed computing methods for the risk measure SR. We review more details of Dunkel and Weber (2010) in the next section.   

Since Problem (\ref{estioce}) is a standard stochastic optimization problem, we
adopt the natural approach SAA to solve it. 
Shapiro et al. (2014) discussed the computing of risk measures based on their stochastic optimization theory built. Especially, they discussed the evaluation of the CVaR and built asymptotic properties of the SAA estimators. The SAA estimation of CVaR was also studied in Trindade et al.
(2007).  In this section, we follow Trindade et al.
(2007) and Shapiro et al. (2014)  and study how to generalize the SAA approach to
estimate the OCE risk measures. Suppose that $\xi_1,
\xi_2, \cdots, \xi_n$ are independent and identically distributed (i.i.d.) observations of $\xi$. The SAA
approach suggests using the following optimization problem
\begin{eqnarray}\label{estiocesaa}
&\displaystyle\mathop{\mbox{minimize}}_{t\in T} & t +
\frac{1}{n}\sum^n_{i=1}l(\xi_i-t)
\end{eqnarray}
to approximate (\ref{estioce}). Let $v_n$ and $v^*$ denote the optimal
objective values of Problems (\ref{estiocesaa}) and (\ref{estioce}) respectively, and $S_n$ and $S$ denote the corresponding sets of optimal solutions. We use $v_n$ as an estimator of $v^*$.

In what follows, we study some statistical properties of the
estimators $v_n$ and $S_n$, based on the theory of SAA. Let $m^2(t) =
\E[l^2(\xi-t)]$. We make the following assumption.
\begin{assumption}\label{ass:varianceoce}
$\sup_{t\in T}m^2(t) < +\infty$.
\end{assumption}

Let $D\left(S_n,S\right)=\sup_{x\in S_n} d(x, S)$ where $d(x, S)=\inf_{y\in S} \|x-y\|$. Note that $D\left(S_n,S\right)$ is the deviation of $S_n$ to $S$. The following result builds the consistency of the estimation, and is straightforward from Theorem 5.4 in Shapiro et al. (2014). 
\begin{proposition}\label{prop:oce-con} 
Suppose that Assumption \ref{ass:varianceoce}
is satisfied. Then
$v_n \rightarrow v^*$ w.p.1, and $\mathbb{D}\left(S_n,S\right)\to 0$ w.p.1. If, moreover, $S$ is a singleton $\left\{t^*\right\}$, then for any $t_n \in S_n$, $t_n\to t^*$ w.p.1.  
\end{proposition}

The consistency guarantees that the estimator $v_n$ is valid in a statistical
sense, and the estimation error diminishes as the computational
effort increases. Now we study the convergence rate of $v_n$. The theory is mainly based on
Shapiro et al. (2014). Define the random variable
$$
Y(t) := \lim_{n\rightarrow\infty}
\sqrt{n}\left[\left(t+\frac{1}{n}\sum^n_{i=1}l(\xi_i-t)\right)-\left(t+\E[l(\xi-t)]\right)\right].
$$
Then it follows from the central limit theorem (CLT) (Durrett 2019) that $Y(t)$ has a
normal distribution $N(0, \sigma^2(t))$ where $\sigma^2(t) =
\hbox{Var}[l(\xi-t)]$. Recall that $S$ is the set of optimal solutions of
Problem (\ref{estioce}). Rockafellar and Royset (2015) showed that
$S$ is a compact interval. Suppose the ending points of $S$ are $t_l$ and $t_u$, i.e., $S=[t_l,t_u]$. We make the following assumption.   
\begin{assumption}\label{ass:Lipsoce}
There exists a random variable $K(\xi)$, such that
$\E\left[K(\xi)^2\right]<\infty$ and
\begin{equation*}\label{equ:Lipsoce}
\left|l(\xi-t_1)-l(\xi-t_2)\right|\le K(\xi)\left|t_1-t_2\right|
\end{equation*}
for all $t_1, t_2 \in T$ and a.e. $\xi$.
\end{assumption}
Assumption \ref{ass:Lipsoce} is a Lipschitz condition for the loss function which is widely adopted in stochastic optimization. Based on Theorem 5.7 of Shapiro et al. (2014), we have the following
proposition.
\begin{proposition}\label{thm:oce-rate}
Suppose that Assumptions \ref{ass:varianceoce} and \ref{ass:Lipsoce}
are satisfied. Then
\begin{eqnarray*}
\sqrt{n}(v_n - v^*) \Rightarrow \inf_{t\in S=[t_l,t_u]} Y(t).
\end{eqnarray*}
If, moreover, $S$ is a singleton $\left\{t^*\right\}$, then
\begin{eqnarray*}
\sqrt{n}(v_n - v^*) \Rightarrow N(0, \sigma^2(t^*)).
\end{eqnarray*}
\end{proposition}

Proposition \ref{thm:oce-rate} shows that the convergence rate of the
SAA estimator for OCE risk measure is in the order of
$n^{-1/2}$. It includes two cases for the asymptotic
distribution. In the first case, the set of optimal solutions is a
non-degenerate compact interval and the asymptotic
distribution is an infimum of a set of normal distributions. An
instance for this case is the CVaR. It is well known that the set of
optimal solutions of Problem (\ref{estioce}) for the CVaR setting
may be a compact interval with the corresponding VaR being an end
point of the interval. The results were analyzed in Trindade et al.
(2007). In the second case, the optimal solution is unique and the
asymptotic distribution is simply a normal distribution. This is the
result of asymptotic normality, based on which the approximate
confidence intervals may be built accordingly. In many cases, the
uniqueness of the solution for (\ref{estioce}) can be guaranteed.
For instance, if the objective function in (\ref{estioce}) is
strictly convex, the optimal solution is unique. For the CVaR, if
the random loss has positive density in a neighborhood of the VaR,
then the solution is unique. For this scenario, Proposition 
\ref{thm:oce-rate} recovers the asymptotic normality of the CVaR
estimator built in Trindade et al. (2007).

Suppose the uniqueness of the solution is guaranteed. To construct confidence intervals for the
OCE risk measure, we need to estimate $\sigma^2(t^*)$ where
$t^*$ is the unique solution. Given an i.i.d. sample $\xi_1, \xi_2,
\cdots, \xi_n$ from $\xi$, we solve Problem  (\ref{estiocesaa}) to obtain an optimal solution
$t_n$. We then propose the following estimator for $\sigma^2(t^*)$
\begin{equation*}
\label{eq:var-est} \sigma^2_n := \frac{1}{n}\sum^n_{i=1}l^2(\xi_i-
t_n) -\left(\frac{1}{n}\sum^n_{i=1}l(\xi_i- t_n)\right)^2.
\end{equation*}
We can show the
following consistency of $\sigma^2_n$.
\begin{proposition}\label{prop:est-var}
Suppose that Assumptions \ref{ass:varianceoce} and \ref{ass:Lipsoce}
are satisfied. Then $\sigma^2_n\rightarrow \sigma^2(t^*)$ w.p.1.
\end{proposition}
\noindent\textbf{Proof.} 
By strong
law of large numbers (Durrett 2019), 
$\frac{1}{n}\sum^n_{i=1} l^2\left(\xi_i - t^*\right)\to
\E\left[l^2\left(\xi-t^*\right)\right]$ w.p.1. By Proposition \ref{prop:oce-con}, $t_n\to t^*$ w.p.1. Define
$a_n(t)=\frac{1}{n}\sum^n_{i=1} l^2\left(\xi_i - t\right)$. From Theorem 7.53 of Shapiro et al. (2014), we have $a_n(t)$
converges uniformly over $T$ to
$a(t)=\E\left[l^2\left(\xi-t\right)\right]$ w.p.1.

Let $w$ denote a sample path of $\left\{\xi_1,\xi_2,\cdots\right\}$, and $\Omega$ denote the
set of sample paths for which $t_n\to t^*$, $\frac{1}{n}\sum^n_{i=1}
l^2\left(\xi_i - t^*\right)\to
\E\left[l^2\left(\xi-t^*\right)\right]$ and $a_n(t)\to a(t)$
uniformly. Then the probability of $\Omega$ is $1$. We show that for any $\omega\in \Omega$, $a_n(t_n)\to a(t^*)$.
Fixing $\omega\in \Omega$, we obtain a deterministic sequence, for
which the following inequality holds
\begin{equation}\label{equ:consi:sigma}
\left|a_n(t_n)-a(t^*)\right|\le
\left|a_n(t_n)-a(t_n)\right|+\left|a(t_n)-a(t^*)\right|.
\end{equation}
Since $a_n(t)$ converges uniformly over $T$ to $a(t)$, we have $\left|a_n(t_n)-a(t_n)\right|\to 0$ as $n\to
\infty$. Because $a(t)$ is continuous in $t$,we have $\left|a(t_n)-a(t^*)\right|\to 0$ as
$n\to \infty$. Thus $a_n(t_n)\to a(t^*)$. This implies $\frac{1}{n}\sum^n_{i=1}
l^2\left(\xi_i - t_n\right)\to
\E\left[l^2\left(\xi-t^*\right)\right]$ w.p.1.

Similarly, we can show that $\frac{1}{n}\sum^n_{i=1}l(\xi_i- t_n)\to
\E\left[l\left(\xi-t^*\right)\right]$ w.p.1, and thus we have
\begin{equation*}  
\left(\frac{1}{n}\sum^n_{i=1}l(\xi_i- t_n)\right)^2\to
\left(\E\left[l\left(\xi-t^*\right)\right]\right)^2
\end{equation*}
 w.p.1. Recalling the definition of variance, we have $\sigma^2_n\rightarrow \sigma^2(t^*)$ w.p.1. \hfill$\Box$

Proposition \ref{prop:est-var} guarantees that we can use
$\sigma^2_n$ to estimate $\sigma^2(t^*)$. For a confidence level
$\alpha$, we can build the following (approximate) $1-\alpha$
confident interval for $v^*$
$$
\left[v_n - z_{1-\alpha/2}\sigma_n/\sqrt{n}, v_n +
z_{1-\alpha/2}\sigma_n/\sqrt{n}\right]
$$
where $z_{1-\alpha/2}$ is the $1-\alpha/2$ quantile of the standard
normal distribution.

Suppose that $S$ is a singleton $\left\{t^*\right\}$. We also build a CLT for the optimal solution. As discussed above, Hamm et al. (2013) used some stochastic approximation procedure to estimate $t^*$. They built a CLT for the estimator. In this paper, we consider the SAA. Let $t_n$ denote an optimal solution of Problem (\ref{estiocesaa}) and use it as the estimator of $t^*$. We build the following proposition. 
\begin{proposition}\label{pro:CLTOCEsolution}
Suppose that Assumptions \ref{ass:varianceoce} and \ref{ass:Lipsoce} are satisfied. Suppose that $l(\xi-t)$ is continuously differentiable w.p.1 and that $\left.\frac{d}{dt}\E\left[l'\left(\xi-t\right)
\right]\right|_{t=t^*}<0$. Then 
\begin{equation}\label{CLTOCEsolution}
\sqrt{n}\left(t_n-t^*\right)\Rightarrow  N(0, \sigma^2), 
\end{equation}
where 
\begin{equation*}\label{}
\sigma^2=\frac{\mathrm{Var}\left[
l'\left(\xi-t^*\right)\right]}{\left[\left.\frac{d}{dt}\E\left[l'\left(\xi-t\right)
\right]\right|_{t=t^*}\right]^2}.
\end{equation*}
\end{proposition}
\noindent\textbf{Proof.} 
The loss function $l$ is convex. Then $l'$ is nondecreasing. By the assumptions, we have  
\begin{equation*}
 \frac{d}{dt}\E\left[l\left(\xi-t\right)\right]=-\E\left[l'\left(\xi-t\right)
\right].
\end{equation*}
Furthermore, given $\frac{d}{dt}\E\left[l'\left(\xi-t^*\right)
\right]<0$, it follows from Corollary 3.5 of Huber and Ronchetti (2011) that (\ref{CLTOCEsolution}) holds. This concludes the proof. \hfill$\Box$

In Proposition \ref{pro:CLTOCEsolution}, because $l'$ is nondecreasing, the condition that $\left.\frac{d}{dt}\E\left[l'\left(\xi-t\right)
\right]\right|_{t=t^*}<0$ could be easily satisfied. In the CLT built in (\ref{CLTOCEsolution}), the asymptotic variance involves $\left.\frac{d}{dt}\E\left[l'\left(\xi-t\right)
\right]\right|_{t=t^*}$. If $l'\left(\xi-t\right)$ is differentiable for a.e. $\xi$ and the Lipschitz type condition (similar as Assumption \ref{ass:Lipsoce}) holds for $l'\left(\xi-t\right)$, we can further put the derivative into the expectation and obtain 
\begin{equation}\label{OCECLTvariance2}
\sigma^2=\frac{\mathrm{Var}\left[
l'\left(\xi-t^*\right)\right]}{\left[\E\left[l''\left(\xi-t^*\right)
\right]\right]^2}.
\end{equation} 
To estimate (\ref{OCECLTvariance2}) we can estimate the numerator and the
denominator in (\ref{OCECLTvariance2}) separately. For the two values, we can use the plug-in estimators as for $\sigma^2(t^*)$ in Proposition \ref{prop:est-var}. That is, we first estimate $t^*$ by $t_n$, and then  plug in $t_n$ and estimate $\mathrm{Var}\left[
l'\left(\xi-t_n\right)\right]$ and $\left[\E\left[l''\left(\xi-t_n\right)
\right]\right]^2$ using the estimators of variance and mean. The consistency of the estimation may be built under some conditions by similar argument as that for Proposition \ref{prop:est-var}.

\subsection{Sensitivity Analysis of OCE Risk Measure}

We now consider the sensitivity of the OCE risk measure. Suppose
that the random loss of interest is $\xi(\theta)$ which depends on a
parameter $\theta$. For simplicity of notation let $\rho(\theta)$ denote the OCE risk measure of $\xi(\theta)$. Then by the definition, $\rho(\theta)$ is the optimal value of the following optimization problem
\begin{eqnarray}\label{estioce_theta}
&\displaystyle\mathop{\mbox{minimize}}_{t\in T} & t +\E[l(\xi(\theta)-t)].
\end{eqnarray}
We aim to evaluate the derivative
$d \rho (\theta)/d\theta$. We assume that $\theta$ is one-dimensional and we have a compact set $\Theta\subset \mathbb{R}$ such that $\theta\in \Theta$. Let $\xi'(\theta)$ be the sample path derivative and assume that it exists w.p.1. and can be evaluated readily, see, e.g., Hong and Liu (2009) for a more detailed discussion on this. We make the following assumptions. 

\begin{assumption}\label{ass:ocesen1}
For every $\theta\in \Theta$, Problem (\ref{estioce_theta}) has a unique optimal solution $t^*(\theta)$.  The compact interval $T$ is specified such that $t^*(\theta)$ is an interior point of $T$ for all $\theta\in \Theta$. 
\end{assumption}

\begin{assumption}\label{ass:ocesen2}
For every $(\theta,t) \in \Theta\times T$, $l(\xi(\theta)-t)$ is continuously differentiable at $(\theta,t) \in \Theta\times T$ w.p.1. There exists a random variable $L$, such that $\E\left[L^2\right]<\infty$ and
\begin{equation*}
|l(\xi(\theta)-t)|\le L  
\end{equation*}
for all $(\theta,t)\in \Theta\times T$ and a.e. $\xi$. There exists a random variable $K$, such that
$\E\left[K^2\right]<\infty$ and
\begin{equation*}\label{equ:Lips}
\left|l(\xi(\theta_1)-t_1)-l(\xi(\theta_2)-t_2)\right|\le
K\left\|(\theta_1,t_1)-(\theta_2,t_2)\right\|
\end{equation*}
for all $(\theta_1,t_1), (\theta_2,t_2) \in \Theta\times T$ and a.e. $\xi$.
\end{assumption}

Assumptions \ref{ass:ocesen1} and \ref{ass:ocesen2} guarantee the smoothness of the OCE risk measure. 
The following proposition provides an expression for the derivative of the OCE risk measure. 
\begin{proposition}\label{pro:OCE:sen}
Suppose that Assumptions \ref{ass:ocesen1} and \ref{ass:ocesen2} are satisfied. Then for $\theta$,  
\begin{equation}\label{OCE-derivative-close}
\frac{d\rho(\theta)}{d\theta}=\E\left[l'(\xi(\theta)-
t^*(\theta))\xi'(\theta)\right].
\end{equation}
\end{proposition}
\noindent\textbf{Proof.} 
Since Assumption \ref{ass:ocesen2} is satisfied, we have that for any $t\in T$, $t +\E[l(\xi(\theta)-t)]$ is differentiable in $\theta$, and $\frac{d}{d \theta}\left[t +\E[l(\xi(\theta)-t)]\right]$ is continuous in $(\theta,t)$. Note that $t^*(\theta)$ is the unique optimal solution. Then by the Danskin theorem (Shapiro et al. 2014),  
\begin{equation*}\label{}
\frac{d\rho(\theta)}{d\theta}=\left.\frac{d}{d \theta}\left[t +\E[l(\xi(\theta)-t)]\right]\right|_{t=t^*(\theta)}.
\end{equation*}
Given Assumptions \ref{ass:ocesen1} and \ref{ass:ocesen2}, by Theorem 7.57 of Shapiro et al. (2014), we have 
\begin{equation*}\label{}
\frac{d}{d \theta}\left[t +\E[l(\xi(\theta)-t)]\right]=\E\left[l'(\xi(\theta)-
t)\xi'(\theta)\right].
\end{equation*}
The two equations above imply (\ref{OCE-derivative-close}). The proof is completed. \hfill$\Box$

When the OCE risk measure is the CVaR, $l=\frac{1}{\alpha}[\cdot]^+$ and $t^*(\theta)$ is the VaR of $\xi(\theta)$. We denote $t^*(\theta)=\mbox{VaR}_{1-\alpha}(\xi(\theta))$.  Then 
 \begin{equation}\label{CVaR-derivative-close}
\frac{d\rho(\theta)}{d\theta}=\frac{1}{\alpha}\E\left[1_{\left\{\xi(\theta)-
t^*(\theta)\ge 0\right\}}\xi'(\theta)\right].
\end{equation}
This expression was derived in Hong and Liu (2009). They used a different approach to build the expression but also outlined the approach of using the Danskin theorem. Note that Hong and Liu (2009) made a slightly different set of assumptions to derive (\ref{CVaR-derivative-close}). They assumed that $\xi(\theta)$ has a density in the neighborhood of the VaR of $\xi(\theta)$. This ensures  Problem (\ref{estioce_theta}) has a unique optimal solution $t^*(\theta)$,  in accordance with Assumption \ref{ass:ocesen1}. They also made a Lipschitz condition on $\xi(\theta)$ that could be covered by Assumption \ref{ass:ocesen2}. In Proposition \ref{pro:OCE:sen}, we extend the result for a more general OCE risk measure. 

We can estimate $\frac{d\rho(\theta)}{d\theta}$ based on Proposition \ref{pro:OCE:sen}. We first use a set of i.i.d. sample $\left\{\xi_1,\cdots, \xi_m\right\}$ to estimate $t^*(\theta)$. Denote the estimator as $t_m(\theta) $. We then take another set of i.i.d. sample $\left\{\xi_1,\cdots, \xi_n\right\}$ which is independent of that for computing $t_m(\theta) $. We use
\begin{equation*}
T_n= \frac{1}{n} \sum_{i=1}^n l'(\xi_i(\theta)-t_m(\theta) )\xi_i'(\theta) 
\end{equation*}
 as the estimator for $\frac{d\rho(\theta)}{d\theta}$. We build the following CLT for $T_n$. 
\begin{proposition}\label{pro:derivativeCLT}
Suppose that Assumptions \ref{ass:ocesen1} and \ref{ass:ocesen2} are satisfied, the loss function $l$ is twice continuously differentiable and $\E\left[l''\left(\xi-
t^*(\theta)\right)\xi'(\theta)\right]$ is finite valued, and that $\left.\frac{d}{dt}\E\left[l'\left(\xi(\theta)-t\right)\right]\right|_{t=t^*(\theta)}<0$.  
Suppose that $m/n \to +\infty$ as $n\to +\infty$. Then  
\begin{equation}\label{equ:derivativeCLT}
\sqrt{n}\left(T_n-\frac{d\rho(\theta)}{d\theta}\right)\Rightarrow  N(0, \sigma^2), n\to +\infty, 
\end{equation}   
where $\sigma^2=\hbox{Var}[l'(\xi(\theta)-
t^*(\theta))\xi'(\theta)]$.
\end{proposition}
\noindent\textbf{Proof.} For simplicity of notation, let $\xi=\xi(\theta)$ and $\zeta=\xi'(\theta)$. Note that for each $i$, the Taylor
expansion yields that
\begin{equation*}
l'\left(\xi_i-t_m(\theta)\right)\zeta_i=l'\left(\xi_i-
t^*(\theta)\right)\zeta_i-l''\left(\xi_i-
t^*(\theta)\right)\zeta_i\left[t_m(\theta)-t^*(\theta)\right]+o\left(t_m(\theta)-t^*(\theta)\right).
\end{equation*}
It follows that
\begin{eqnarray}
\sqrt{n}\left(T_n-\frac{d\rho(\theta)}{d\theta}\right) &=& \sqrt{n}
\left[\frac{1}{n}\sum_{i=1}^nl'\left(\xi_i-
t^*(\theta)\right)\zeta_i-\frac{d\rho(\theta)}{d\theta}\right]\nonumber\\
&& -\left[\frac{1}{n}\sum_{i=1}^n l''\left(\xi_i-
t^*(\theta)\right)\zeta_i\right]\sqrt{n}\left[t_m(\theta)-t^*(\theta)\right]\nonumber \\
&& +\sqrt{n} \
o\left(t_m(\theta)-t^*(\theta)\right).\label{equ:CLT1}
\end{eqnarray}
By Assumption \ref{ass:ocesen2}, $\sigma^2=\hbox{Var}[l'(\xi(\theta)-
t^*(\theta))\xi'(\theta)]<+\infty$. By Proposition \ref{pro:OCE:sen} and the CLT,
\begin{equation*}
\sqrt{n} \left[\frac{1}{n}\sum_{i=1}^nl'\left(\xi_i-
t^*(\theta)\right)\zeta_i-\frac{d\rho(\theta)}{d\theta}\right]\Rightarrow N(0, \sigma^2).
\end{equation*}
Because $m/n \to +\infty$ as $n\to +\infty$, it follows from Proposition \ref{pro:CLTOCEsolution} that
$\sqrt{n}\left(t_m(\theta)-t^*(\theta)\right)\to 0$ in probability as $n\to
+\infty$. By strong law of large numbers, $\frac{1}{n}\sum_{i=1}^n
l''\left(\xi_i- t^*(\theta)\right)\zeta_i\to \E\left[l''\left(\xi-
t^*(\theta)\right)\zeta\right]$ w.p.1. Then we have the second term
on the right hand side of (\ref{equ:CLT1}) converges to
$0$ in probability as $n\to +\infty$. Furthermore, $\sqrt{n} \
o\left(t_m(\theta)-t^*(\theta)\right)\to 0$  in probability as $n\to +\infty$.
Therefore, by Slutsky's theorem, (\ref{equ:derivativeCLT}) holds. 
This concludes the proof of the proposition. \hfill$\Box$

In Proposition \ref{pro:derivativeCLT}, we require that $m$ tends to $+\infty$ faster than $n$. Intuitively, we require that the inner estimator is more accurate. Proposition \ref{pro:derivativeCLT} shows that the CLT holds for the estimator $T_n$. We can also build a confidence region for $\frac{d\rho(\theta)}{d\theta}$ based on the asymptotic result.  

Glynn et al. (2020) studied the CLT for estimated functions at estimated points. They considered the setting where the stochastic process is stochastically continuous. By using their results, it is possible to relax some of the conditions and build the CLT for the estimator of $\frac{d\rho(\theta)}{d\theta}$ accordingly.

\subsection{Optimization of OCE Risk Measure} 

In real decision, we often need to optimize the risk measure. We consider the optimization of the OCE risk measure  
\begin{eqnarray}\label{optoce}
&\displaystyle\mathop{\mbox{minimize}}_{x\in \mathcal{X},\   t\in T} & t +\E[l(c(x,\xi)-t)], 
\end{eqnarray}
where $c$ is a function of the decision vector $x$ and the random vector $\xi$ and $ \mathcal{X}$ is a compact set. With a set of i.i.d. sample $\left\{\xi_1,\cdots, \xi_n\right\}$, the SAA of Problem (\ref{optoce}) is 
\begin{eqnarray*}\label{optoceSAA}
&\displaystyle\mathop{\mbox{minimize}}_{ x\in \mathcal{X},\ t\in T } & \frac{1}{n} \sum_{i=1}^n\left\{t +l(c(x,\xi_i)-t)\right\}. 
\end{eqnarray*}
Define the random variable
$$
Y(x,t) := \lim_{n\rightarrow\infty}
\sqrt{n}\left[\left(t+\frac{1}{n}\sum^n_{i=1}l(c(x,\xi_i)-t)\right)-\left(t+\E[l(c(x,\xi)-t)]\right)\right].
$$
Let $m^2(x,t) =
\E[l^2(c(x,\xi)-t)]$ and $\sigma^2(x,t) =
\hbox{Var}[l(c(x,\xi)-t)]$. Similar to Assumption \ref{ass:varianceoce}, we assume that $\sup_{x\in \mathcal{X}, t\in T}m^2(x,t) < +\infty$. We further make the following assumption. 
\begin{assumption}\label{ass:Lipsoceopt}
There exists a random variable $K(\xi)$, such that
$\E\left[K(\xi)^2\right]<\infty$ and
\begin{equation*}\label{equ:Lipsoce}
\left|l(c(x_1,\xi)-t_1)-l(c(x_2,\xi)-t_2)\right|\le K(\xi)\left\|(x_1,t_1)-(x_2,t_2)\right\|
\end{equation*}
for all $x_1,x_2\in \mathcal{X}$, $t_1, t_2 \in T$ and a.e. $\xi$.
\end{assumption}

Let $S$ be the set of optimal solutions of
Problem (\ref{optoce}). The following result is a direct application of the asymptotic results of Shapiro et al. (2014). 
\begin{proposition}
Suppose that Assumption \ref{ass:Lipsoceopt}
is satisfied. Then
\begin{eqnarray*}
\sqrt{n}(v_n - v^*) \Rightarrow \inf_{(x,t)\in S} Y(x, t).
\end{eqnarray*}
If, moreover, $S$ is a singleton $\left\{(x^*,t^*)\right\}$, then
\begin{eqnarray*}
\sqrt{n}(v_n - v^*) \Rightarrow N(0, \sigma^2(x^*,t^*)).
\end{eqnarray*}
\end{proposition}
From the analysis above, we can see that for both estimation and optimization, the typical $n^{-1/2}$ rate of convergence could be achieved.

\subsection{Smooth Conditional Value-at-risk} 

We have discussed the computing of the OCE risk measure and built various theoretical properties. As can be seen, these properties typically require certain technical conditions. An important condition is an essential smoothness of the risk measure. The CVaR is the most widely used OCE risk measure. However, from a computation perspective, the CVaR may sometimes face the issue of smoothness sufficiency. When optimizing a CVaR function, the smoothing techniques are often used to conduct the optimization. To address the issue, we can directly define alternative smooth risk measures. By defining the loss function $l(z)=\frac{1}{\alpha}\left([z]^+\right)^p, p\ge 1$, we can generate a convex risk measure via  
\begin{eqnarray}\label{SCVaR}
&\displaystyle\mathop{\mbox{minimize}}_{t\in \mathbb{R}} & t +\frac{1}{\alpha}\E[\left([\xi-t]^+\right)^p], p\ge 1. 
\end{eqnarray}
When $p=1$, this risk measure is the CVaR. When $p>1$, the loss function $l$ is continuously differentiable. Thus, by selecting some $p>1$, we may avoid the difficulty of nonsmoothness. The risk measure defined by (\ref{SCVaR}) shares similar structure as the higher-moment coherent risk measures proposed in Krokhmal (2007), which  is given by 
\begin{eqnarray*}\label{SCVaR_Co}
&\displaystyle\mathop{\mbox{minimize}}_{t\in \mathbb{R}} & t +\frac{1}{\alpha}\left[\E[\left([\xi-t]^+\right)^p]\right]^{1/p}, p\ge 1. 
\end{eqnarray*}

We can also define more smooth risk measures. Alexander et al. (2006) studied the optimization of CVaR for the portfolio of derivatives. They focused on smoothing the CVaR function by using a quadratic function to smooth the piecewise linear function. Following their smooth function, we can define the loss function 
 \begin{equation}\label{quadratic_smooth} 
l(z)=\left\{\begin{array}{cc} \frac{z}{\alpha} & z\ge \epsilon \\[6pt]  
\frac{1}{\alpha}\left[\frac{z^2}{4\epsilon}+\frac{1}{2} z +\frac{1}{4}\epsilon\right] & -\epsilon\le z\le \epsilon \\[6pt]  
0 & {\mbox{otherwise}}. 
\end{array}\right.
\end{equation} 
The convex risk measure induced by (\ref{quadratic_smooth}) becomes smoother. Note that when $\epsilon=0$, it becomes the CVaR.

\section{Utility-Based Shortfall Risk}\label{sec:SR}

In the previous section, we discussed the risk measure that can be defined as the optimal value of a stochastic optimization problem that optimizes an expectation function. In this section, we consider some risk measures that can be defined (written) as a stochastic optimization problem with expectation constraint.  Our focus is on the utility-based shortfall risk (SR).  

SR is a widely used convex risk measure, and was introduced in F\"{o}llmer and Schied (2002). Following Dunkel and Weber (2010), let $l$ be a convex loss function and $\lambda$ be a prespecified value. Define 
\begin{equation*}
\mathscr{A}:=\left\{X\in L^{\infty}: \E\left[l(-X)\right]\le
\lambda\right\}
\end{equation*}
where $L^{\infty}$ denotes the set of bounded random variables. 
The SR of a financial position $X$ corresponding
to $l$ and $\lambda$ is then induced by the acceptance set
$\mathscr{A}$ and is defined as
\begin{equation*}
\mbox{SR}_{l,\lambda}(X)=\inf\left\{t\in \mathbb{R}: t+X \in
\mathscr{A}\right\}.
\end{equation*}   
For ease of analysis, we consider the loss position $\xi=-X$ and still denote $\mbox{SR}_{l,\lambda}(\xi)$ as the SR of $\xi$. Let 
\begin{equation*}
g(t):=\E\left[l\left(\xi-t\right)\right]-\lambda. 
\end{equation*}
We make the following assumption. Note that Hu and Zhang (2018) provided justification for some of the assumptions made in this section.  
\begin{assumption}\label{ass:slater}
There exist $t_l,t_u\in \mathbb{R}$, such that $g(t_u)< \lambda$ and
$g(t_l)> \lambda$.
\end{assumption}

Dunkel and Weber (2010) studied using the stochastic root finding approach to estimate the SR. They treated SR as the root of the following equation 
\begin{equation}\label{equ:SRroot}
g(t)=\E\left[l\left(\xi-t\right)\right]-\lambda=0. 
\end{equation}
They constructed some measurable function $\hat Y_t: [0,1]\to \mathbb{R}$ such that $g(t)=\E\left[\hat Y_t(U)\right]$ where $U$ is a random number. A conventional construction is to use the inverse cumulative distribution function of $\xi$. Dunkel and Weber (2010) used importance sampling to construct $\hat Y_t$ to improve the speed of the algorithm. For simplicity, let $Y_n=\hat Y_{t_n}(U_n)$ where $U_n$ is independently drawn from $U$. They proposed the following Robbins-Monro stochastic approximation procedure to search the root $t^*$ of Equation (\ref{equ:SRroot})
\begin{equation*}\label{equ:SA}
t_{n+1}=\Pi \left[t_n+\frac{c}{n^{\gamma}}\cdot Y_n\right], 
\end{equation*}
where $c$ and $\gamma$ are parameters and $\Pi$ represents the projection of $t$ onto $[t_l,t_u]$. They built the following asymptotic results for the estimator $t_n$ under certain conditions 
\begin{equation*}
\begin{array}{cc}  \sqrt{n}\left(t_n-t^*\right)\Rightarrow N(0, \sigma_1^2),\quad \mbox{if} \ \gamma=1;  \\
\\ 
\sqrt{n^\gamma}\left(t_n-t^*\right)\Rightarrow N(0, \sigma_2^2),\quad \mbox{if}\  \gamma\in (\frac{1}{2},1). 
\end{array}
\end{equation*}
where 
\begin{equation*}
\sigma_1^2=\frac{-c^2 \mathrm{Var}\left[
l\left(\xi-t^*\right)\right]}{-2c\E\left[l'\left(\xi-t^*\right)
\right]+1}, \quad \sigma_2^2= \frac{-c \mathrm{Var}\left[
l\left(\xi-t^*\right)\right]}{-2\E\left[l'\left(\xi-t^*\right)
\right]+1}. 
\end{equation*}
Dunkel and Weber (2010) also considered a Polyak-Ruppert averaging procedure to search the root. They built a slightly different convergence result for the estimator.     

Hu and Zhang (2018) studied the computation of SR based on SAA. In this section, we review the approaches of Hu and Zhang (2018). Following the formulation, the SR can be treated as the optimal value of the following optimization problem
\begin{eqnarray}\label{estisr}
&\displaystyle\mathop{\mbox{minimize}}_{t\in T} & t\\
&\mbox{ subject to} & \E\left[l\left(\xi-t\right)\right]\le
\lambda,\nonumber
\end{eqnarray}
where $T$ is a well specified compact interval such that the optimal solution of Problem (\ref{estisr}) is an interior point of $T$, e.g., $T=\left[t_l,t_u\right]$.  

Suppose that we obtain a set of i.i.d. sample $\left\{\xi_1,\cdots, \xi_n\right\}$ from the distribution of $\xi$.  Then, the SAA of Problem (\ref{estisr}) takes the following expression
\begin{eqnarray}\label{estisrsaa}
&\displaystyle\mathop{\mbox{minimize}}_{t\in T} & t\\
&\mbox{ subject to} & \frac{1}{n}\sum_{j=1}^n
l\left(\xi_j-t\right)\le \lambda.\nonumber
\end{eqnarray}
Let $t_n$ and $t^*$ denote the optimal objective values of Problems (\ref{estisrsaa}) and (\ref{estisr}) respectively. Hu and Zhang (2018) proposed to use $t_n$ to estimate $t^*$, i.e., the SR. Let $m^2(t)=\E\left[l^2\left(\xi-t\right)\right]$. Similar to Assumption \ref{ass:varianceoce}, it is assumed that $\sup_{t\in T} m^2(t)<+\infty$. Let $D_l$ denote the set of
points at which $l$ is continuously differentiable. They made the following assumption and built the following result. 

\begin{assumption}\label{ass:Lip}
For any $t\in T$,
\begin{equation*}\label{equ:ass:Lip}
\Pr\left\{\xi- t\in D_l\right\}=1,
\end{equation*}
and there exists a random variable $K(\xi)$, such that
$\E\left[K(\xi)^2\right]<\infty$ and for all $t_1, t_2 \in T$ and
a.e. $\xi$,
\begin{equation*}\label{equ:Lip}
\left|l(\xi-t_1)-l(\xi-t_2)\right|\le K(\xi)\left|t_1-t_2\right|.
\end{equation*}
Moreover, $\Pr\left\{\xi: l'(\xi-t^*)>0\right\}>0$.
\end{assumption}

\begin{theorem}[Hu and Zhang (2018)]\label{thm:estisr:norm}
Suppose that Assumptions \ref{ass:slater} to \ref{ass:Lip} are
satisfied. Then 
\begin{equation*}
\sqrt{n}\left(t_n-t^*\right)\Rightarrow N(0, \sigma^2),
\end{equation*}
where
\begin{equation*}\label{normvariance}
\sigma^2=\frac{\mathrm{Var}\left[
l\left(\xi-t^*\right)\right]}{\left[\E\left[l'\left(\xi-t^*\right)
\right]\right]^2}=\frac{\E\left[l^2\left(\xi-t^*\right)\right]-\lambda^2}{\left[\E\left[l'\left(\xi-t^*\right)
\right]\right]^2}.
\end{equation*}
\end{theorem}
From Theorem \ref{thm:estisr:norm}, it can be seen that the convergence rate of $t_n$ is still  in the order of $n^{-1/2}$. Hu and Zhang (2018) also discussed some interesting connection of the results in Theorem \ref{thm:estisr:norm} with the asymptotic results that are built in Dunkel and Weber (2010). They showed that compared with the stochastic root finding approach without importance sampling in Dunkel and Weber (2010), the SAA result in Theorem \ref{thm:estisr:norm} achieves the best convergence rate or the smallest asymptotic variance.  

For the sensitivity of the SR, we introduce $\xi(\theta), \theta\in \Theta$, and let $t^*(\theta)$ denote the SR of $\xi(\theta)$.  Hu and Zhang (2018) made the following assumptions and derived the expression of the derivative of the SR.  

\begin{assumption}\label{ass:slatersen}
For each $\theta\in\Theta$, there exist $t_l,t_u\in \mathbb{R}$, satisfying
$\E\left[l\left(\xi(\theta)-t_u\right)\right]< \lambda$ and
$\E\left[l\left(\xi(\theta)-t_l\right)\right]> \lambda$.
\end{assumption}

\begin{assumption}\label{ass:smoothsen}
There exists a compact interval $T$ of $t$, such that for any $\theta\in\Theta$, $t^*(\theta)$ is an interior point of $T$, and for any $t\in T$,
\begin{equation*}\label{equ:ass:Lipsen}
\Pr\left\{\xi(\theta)- t\in D_l\right\}=1.
\end{equation*}
Moreover, $\Pr\left\{\xi(\theta):
l'(\xi(\theta)-t^*(\theta))>0\right\}>0$.
\end{assumption}

\begin{assumption}\label{ass:Lips}
For any
$(\theta,t)\in\Theta\times T$, $\partial l(\xi(\theta)- t)/\partial
\theta$ and $\partial l(\xi(\theta)- t)/\partial t$ exist w.p.1 and
\begin{equation*}\label{equ:chainrule}
\frac{\partial l(\xi(\theta)- t)}{\partial \theta}=l'(\xi(\theta)-
t)\xi'(\theta),\ \frac{\partial l(\xi(\theta)- t)}{\partial
t}=-l'(\xi(\theta)- t).
\end{equation*}
Furthermore, there exists a random variable $K$, such that
$\E\left[K\right]<\infty$ and
\begin{equation*}\label{equ:Lips}
\left|l(\xi(\theta_1)-t_1)-l(\xi(\theta_2)-t_2)\right|\le
K\left\|(\theta_1,t_1)-(\theta_2,t_2)\right\|
\end{equation*}
for all $(\theta_1,t_1), (\theta_2,t_2) \in \Theta\times T$.
\end{assumption}

\begin{theorem}[Hu and Zhang (2018)]
Suppose that Assumptions \ref{ass:slatersen} to \ref{ass:Lips} are satisfied. Then, 
\begin{equation*}\label{derivative-close}
\frac{dt^*(\theta)}{d\theta}=\frac{\E\left[l'(\xi(\theta)-
t^*(\theta))\xi'(\theta)\right]}{\E\left[l'(\xi(\theta)-
t^*(\theta))\right]}.
\end{equation*}
\end{theorem}
Hu and Zhang (2018) constructed an estimator for $\frac{dt^*(\theta)}{d\theta}$ and built the consistency of the estimator. They also discussed how to build a confidence interval for $\frac{dt^*(\theta)}{d\theta}$. 

We further analyze the optimization of SR. Consider the following optimization problem
\begin{eqnarray}
&\displaystyle\mathop{\mbox{minimize}}_{x\in \mathcal{X}} & \mbox{SR}_{l,\lambda}(c(x,\xi)), \nonumber
\end{eqnarray}
where $\mathcal{X}$ is a convex compact set. The problem can be transformed as 
\begin{eqnarray}\label{optsr}
&\displaystyle\mathop{\mbox{minimize}}_{x\in \mathcal{X}, \  t\in T} & t\\
&\mbox{ subject to} & \E\left[l\left(c(x,\xi)-t\right)\right]\le
\lambda, \nonumber
\end{eqnarray}
where $T$ is a well specified compact interval such that the optimal $t$ is an interior point of $T$. Given the i.i.d. sample $\left\{\xi_1,\cdots, \xi_n\right\}$, we construct the following SAA of Problem (\ref{optsr})
 \begin{eqnarray}\label{optsrSAA}
&\displaystyle\mathop{\mbox{minimize}}_{x\in \mathcal{X}, \ t\in T} & t\\
&\mbox{ subject to} & \frac{1}{n}\sum_{i=1}^n\left[l\left(c(x,\xi_i)-t\right)\right]\le
\lambda. \nonumber
\end{eqnarray}
Let $t_n$ and $t^*$ denote the optimal objective
values of Problems (\ref{optsrSAA}) and (\ref{optsr}) respectively. We build a CLT for $t_n$ based on the theory of SAA. This was not explored in Hu and Zhang (2018). We make the following assumptions.  

\begin{assumption}\label{ass:convexSRopt}
The function $c(x,\xi)$ is convex in $x$ for a.e. $\xi$. 
\end{assumption}

\begin{assumption}\label{ass:Slater}
There exists $x\in {\mbox int} \ \mathcal{X}$ and $t\in {\mbox int} \ T$ such that $\E\left[l\left(c(x,\xi)-t\right)\right]< \lambda$. 
\end{assumption}

Assumption \ref{ass:Slater} is a Slater condition which is a standard constraint qualification in optimization.  

\begin{assumption}\label{ass:LipsSRopt}
The function $l(c(x,\xi)-t)$ is differentiable with respect to $(x,t)$ w.p.1. There exists a random variable $K(\xi)$, such that
$\E\left[K(\xi)^2\right]<\infty$ and
\begin{equation*}\label{equ:Lipsoce}
\left|l(c(x_1,\xi)-t_1)-l(c(x_2,\xi)-t_2)\right|\le K(\xi)\left\|(x_1, t_1)-(x_2, t_2)\right\|
\end{equation*}
for all $x_1,x_2\in \mathcal{X}$, $t_1, t_2 \in T$ and a.e. $\xi$.
\end{assumption}

\begin{assumption}\label{ass:uniqueSRopt}
Problem (\ref{optsr}) has a unique optimal solution $(x^*,t^*)$, and
\begin{equation*}
\Pr\left\{\xi: l'(c(x^*,\xi)-t^*)>0\right\}>0.
\end{equation*}
\end{assumption}

We build the following theorem. 
\begin{theorem}\label{thm:optsr:norm}
Suppose that Assumptions \ref{ass:convexSRopt} to \ref{ass:uniqueSRopt} are
satisfied. Then
\begin{equation*}
\sqrt{n}\left(t_n-t^*\right)\Rightarrow N(0, \sigma^2),
\end{equation*}
where
\begin{equation}\label{normvariance}
\sigma^2=\frac{\mathrm{Var}\left[
l\left(c(x^*,\xi)-t^*\right)\right]}{\left[\E\left[l'\left(c(x^*,\xi)-t^*\right)
\right]\right]^2}=\frac{\E\left[l^2\left(c(x^*,\xi)-t^*\right)\right]-\lambda^2}{\left[\E\left[l'\left(c(x^*,\xi)-t^*\right)
\right]\right]^2}.
\end{equation}
\end{theorem}

\noindent\textbf{Proof.} We follow a similar argument as for the above Theorem \ref{thm:estisr:norm} which was provided by Hu and Zhang (2018).       
Consider the Lagrangian duality of
Problem (\ref{optsr}). The associated Lagrangian is
\begin{equation*}
L(x,t,\mu)=t+\mu\left[\E\left[l\left(c(x,\xi)-t\right)\right]-\lambda\right].
\end{equation*}
Then the Lagrangian dual of Problem  (\ref{optsr}) is 
\begin{equation*}
\displaystyle\mathop{\mbox{maximize}}_{\mu\in \mathbb{R}^+} \
\displaystyle\mathop{\mbox{minimize}}_{x\in \mathcal{X}, \ t\in T} \ L(x,t,\mu).
\end{equation*}

Let $S$ be the set of optimal solutions of Problem (\ref{optsr}) and
$\Lambda$ be the set of associated multipliers. By assumption $S=\{(x^*, t^*)\}$ is a singleton. By Assumption \ref{ass:LipsSRopt}, for each $\mu\ge 0$,
$L(x,t,\mu)$ is differentiable w.r.t. $t$. Note that $t^*$ is an interior point of $T$. We have that the optimal
multiplier $\mu^*$ associated with $(x^*,t^*)$ satisfies
\begin{equation*}
\left.\frac{d L(x,t,\mu^*)}{d t}\right|_{x=x^*, t=t^*}=0.
\end{equation*}
This implies that
\begin{equation*}
1-\mu^*\E\left[l'\left(c(x^*,\xi)-t^*\right)\right]=0,
\end{equation*}
where we can put the derivative into the expectation due to Assumption \ref{ass:LipsSRopt}. From Assumption \ref{ass:uniqueSRopt} and the nondeceasing of $l$, we have
$\E\left[l'\left(c(x^*,\xi)-t^*\right)\right]>0$. Therefore, $\mu^*$ is
uniquely determined by
$\mu^*=1/\E\left[l'\left(c(x^*,\xi)-t^*\right)\right]$. This means
$\Lambda=\{\mu^*\}$ is also a singleton. It follows from Theorem
5.11 of Shapiro et al. (2014) that
$\sqrt{n}\left(t_n-t^*\right)\Rightarrow N(0,\sigma^2)$, where
\begin{eqnarray*}
\sigma^2&:=& \mathrm{Var}\left[\mu^*l\left(c(x^*,\xi)-t^*\right)\right]={\mu^*}^2
\mathrm{Var}\left[l\left(c(x^*,\xi)-t^*\right)\right]\\
&=& \frac{\E\left[l^2\left(
c(x^*,\xi)-t^*\right)\right]-\left[\E\left[l\left(c(x^*,\xi)-t^*\right)\right]
\right]^2}{\left[\E\left[l'\left(c(x^*,\xi)-t^*\right) \right]\right]^2}.
\end{eqnarray*}
Then (\ref{normvariance}) follows from 
$\E\left[l\left(c(x^*,\xi)-t^*\right)\right]=\lambda$. The equation holds since that $\E\left[l\left(c(x^*,\xi)-t\right)\right]$ is convex and thus continuous in $t$. If $\E\left[l\left(c(x^*,\xi)-t^*\right)\right]<\lambda$, we can find $t<t^*$ such that $\E\left[l\left(c(x^*,\xi)-t\right)\right]\le \lambda$. This contradicts with the optimality of $t^*$. The proof is completed. \hfill$\Box$

In Theorem \ref{thm:optsr:norm}, we make the assumption that Problem (\ref{optsr}) has a unique optimal solution $(x^*,t^*)$. The uniqueness of $t^*$ can always be guaranteed because it is the optimal value of Problem (\ref{optsr}). A sufficient condition for the uniqueness of $x$ is that $\E\left[l\left(c(x,\xi)-t\right)\right]$ is strictly convex in $(x, t)$.     
As can be seen from the analysis above, the convergence of optimal value is in the order of
$n^{-1/2}$, which is the conventional convergence rate for a Monte Carlo sampling method.   

The studies of Dunkel and Weber (2010) and Hu and Zhang (2018) built the asymptotic results for the SR. In contrast to the asymptotic analysis, Hegde et al. (2024) studied online estimation and optimization for SR. For the estimation of SR, they treated it as a stochastic root finding problem and used a stochastic approximation scheme similar to Dunkel and Weber (2010) to search the root. Their algorithm recursion is as follows 
\begin{equation*}\label{equ:SA2}
t_{n+1}=\Pi \left[t_n+a_{n+1} \hat g (t_n)\right], 
\end{equation*}
where $\hat g(t)=l(\xi_{n+1}-t)-\lambda$ and $a_{n+1}$ is the step size. They considered a finite sample performance and established some interesting finite sample error bounds under a set of conditions. For the optimization of SR, they proposed a stochastic gradient algorithm and also built finite sample error bounds for their algorithm.

\section{Chance Constrained Programs}\label{sec:CCP}

We have discussed the computing of two important classes of risk measures. In this section, we discuss another widely used risk management model, the chance constrained program (CCP). The CCP is an important stochastic optimization model. From a risk management perspective, it uses the probability, or equivalently, the VaR to constrain the risk. There have been extensive studies on this problem in recent years. In this section, we review some of the advancements of this topic. Consider a conventional CCP with the following expression
\begin{eqnarray}\label{eq:jccp}
&\displaystyle\mathop{\mbox{minimize}}_{x\in X} & h(x)\\
&\mbox{ subject to}& {\Pr}_{\sim P_*} \left\{c_1(x,\xi)\le 0,\cdots,c_m(x,\xi)\le 0\right\}\ge
1-\alpha.\label{eq:jcc}
\end{eqnarray}
In Problem (\ref{eq:jccp}), $x$ is the decision vector, belonging to a set $X$, $\xi\in \Xi$ is a random vector modeling the randomness of the problem, $h(x)$ is a real valued function of $x$, and $c_j(x,\xi), j=1,\cdots,m$ are real valued functions of $x$ and $\xi$ which we call cost functions. Constraint (\ref{eq:jcc}) requires that the $m$ constraints within the probability function are satisfied with at least probability $1-\alpha$, which is a natural relaxation of the requirement that the $m$ constraints are satisfied for all $\xi\in \Xi$. Problem (\ref{eq:jccp}) is often classified as joint CCP (JCCP) if $m\ge 2$ and single CCP (SCCP) if $m=1$. Sometimes, it is beneficial to reformulate Problem (\ref{eq:jccp}) as the following SCCP  
 \begin{eqnarray}\label{eq:ccp}
&\displaystyle\mathop{\mbox{minimize}}_{x\in X} & h(x)\\
&\mbox{ subject to}& {\Pr}_{\sim P_*} \left\{c(x,\xi)\le 0\right\}\ge
1-\alpha,\nonumber\label{eq:cc}
\end{eqnarray}
by defining $c(x,\xi):=\max\left\{c_1(x,\xi),\cdots,c_m(x,\xi)\right\}$. Based on the definition of VaR, Problem (\ref{eq:ccp}) is equivalent to  
\begin{eqnarray*}\label{eq:VaR:lp}
&\displaystyle\mathop{\mbox{minimize}}_{x\in X} & h(x)\\
&\mbox{ subject to} & \mbox{VaR}_{1-\alpha, P_*}(c(x,\xi))\le 0. \label{eq:VaR:lpc}
\end{eqnarray*}

Due to the inherent complexity of the probability function or VaR function, The CCPs are often difficult to handle. Various approaches have been proposed to solve them. The sampling approaches propose to approximate the CCP based on the sample of the random vector. Two sampling approaches are widely studied in the literature. The first one is SAA (e.g., Luedtke and Ahmed 2008). SAA is a conventional stochastic optimization method. As discussed above, in this approach, the sample mean function is used to approximate the true expectation function and the resulting SAA counterpart is solved. Note that when using the SAA, the sample mean function of the probability function becomes a discontinuous step function. Thus, the sample problem usually needs to be reformulated as the integer (or mixed integer) program. Pagnoncelli et al. (2009) investigated the SAA approach for solving the CCP. They built the strong consistency of the SAA optimal value. A second sampling approach is called the scenario approach. It was proposed and studied in Calafiore and Campi (2005, 2006) and De Farias and Van Roy (2004). Hong et al. (2014) reviewed the formulations of these approaches. Küçükyavuz and Jiang (2022) provided a comprehensive review on the sampling based methods for CCP.        

Another natural approach is to treat the CCP as a nonlinear optimization problem and develop corresponding nonlinear optimization procedures to solve it, see, e.g., Henrion and M\"{o}ller (2012) and van Ackooij and Henrion (2014). In contrast to a deterministic problem, to implement a nonlinear optimization algorithm for a CCP, one typically needs to assess the value and the gradient of the probability function in the CCP. Hong (2009) derived gradient formulas of the probability function and the quantile (VaR) function, and developed statistical properties of the quantile estimator. The results could be implemented in a gradient based optimization algorithm to solve the CCP or the VaR constrained program. Hong and Jiang (2017) developed the expressions for both the gradient and the Hessian matrix of a joint probability function. They used the estimators to solve the CCP. Feng and Liu (2016) proposed a change-of-variables (CoV) approach for conditional Monte Carlo and developed an expression which could provide an unbiased estimator for the sensitivities of an expectation function involving indicator function. They implemented the approach to estimate the gradient of a probability function in CCP. They compared their estimator with the conventional  conditional Monte Carlo approach and showed that their approach is more efficient.         

Because the indicator function in the CCP is discontinuous, another smoothing approach proposes to smooth the  indicator function in the probability function. Hong et al. (2011) used a difference-of-convex (DC) function to approximate the indicator function and formulated a DC stochastic program to approximate the CCP. They then proposed a DC algorithm to solve the approximation problem. They built the convergence of the DC approximation as well as the convergence of the DC algorithm. Hu et al. (2013) used a log-sum-exponential function to further smooth the max function in the DC program for a JCCP. Shan et al. (2014) provided a class of smoothing
functions to approximate the CCP. Cui et al. (2022) provided a comprehensive treatment on using the DC functions to approximate the indicator function in the CCP. They also considered the DC cost functions within the probability function and the weighted sum of the probability functions in the CCP, which significantly generalize the model structures. In contrast to dealing with the probability function, Pe\~na-Ordieres et al. (2020) transformed the probabilistic constraint equivalently to a quantile constraint. They adopted the sample approximation for the quantile constraint and developed a smooth approximation approach for the sample problem. 

The chance constraint in the CCP depends on the structure of the cost function and the attribute of the random distribution. Many studies consider some specific structures and distributions in the CCP and make use of these structures and attributes to design solution procedures. An important class of chance constraints has the following expression
\begin{equation}\label{equ:lsep}
{\Pr}_{\sim P_*}\left\{B\xi\le x\right\}\ge 1-\alpha,
\end{equation}
where $B$ is a $m\times k$ matrix and $x$ is the $m$-dimensional decision vector. The CCP with (\ref{equ:lsep}) is a JCCP and many management problems can be formulated as the problem. Henrion and M\"{o}ller (2012) studied this JCCP with the random vector following a multivariate normal distribution. They built a gradient expression of the probability function under this setting and reduced the assessment of the gradient to the assessment of a distribution function of a lower dimensional normal distribution. We introduce their gradient expression below following the introduction of Hu et al. (2022), see the supplemental online materials therein.  

Suppose that in (\ref{equ:lsep}) $P_*$ is the normal distribution $\mathcal {N}(\left.z\right|\mu,\Sigma)$ with mean vector $\mu$ and covariance matrix $\Sigma$. Let $x_i$ denote the $i$-th element of $x$ and  $b_i^\top$ denote the $i$-th row of $B$, and 
\begin{equation}\label{equ:HenrionL1}
S^{(i)}=\Sigma-\frac{1}{b_i^\top\Sigma b_i}\Sigma b_ib_i^\top \Sigma, \ \ w^{(i)}=\mu+\frac{x_i-b_i^\top \mu}{b_i^\top\Sigma b_i}\Sigma b_i,\ \ i=1,\cdots, m.
\end{equation}
From Lemma 3.1 of Henrion and M\"{o}ller (2012), there exists a factorization
\begin{equation}\label{equ:HenrionL2}
S^{(i)}=L^{(i)}{L^{(i)}}^\top,
\end{equation}
where $L^{(i)}$ is a $k \times (k-1)$ matrix with $\mbox{rank}\ L^{(i)} = k-1$. Let
\begin{equation*}
\mathcal{I}(B,x)=\left\{\left.I\subset \left\{1,\cdots,m\right\}\right|\exists z: b_i^\top z = x_i \ (i\in I), b_i^\top z < x_i \ (i\in \left\{1,\cdots,m\right\} \backslash I) \right\}.
\end{equation*}
The system $Bz\le x$ is called nondegenerate if $\mbox{rank} \left\{b_i\right\}_{i\in I}=\# I, \forall I\in \mathcal{I}(B,x)$, where $\#$ denotes the number of elements in a set. We present the following theorem.
\begin{theorem}[Theorem 3.3 of Henrion and M\"{o}ller (2012)]\label{thm:gradient}
Let $x$ be such that $Bz\le x$ is nondegenerate. Suppose that $\Sigma$ is positive definite. Then
\begin{equation*}
\frac{\partial}{\partial x_i} {\Pr}\left\{B\xi\le x\right\}=g_i(x_i) {\Pr}\left\{B^{(i)}L^{(i)}\xi^{(i)}\le x^{(i)}-B^{(i)}w^{(i)}\right\},
\end{equation*}
where $\xi^{(i)}\sim \mathcal {N}(\left.z\right|0,I_{k-1})$, $B^{(i)}$ results from $B$ by deleting row $i$, $x^{(i)}$ results from $x$ by deleting component $i$, $L^{(i)}$ and $w^{(i)}$ are defined in (\ref{equ:HenrionL1})-(\ref{equ:HenrionL2}), and $g_i$ is the one-dimensional normal density with mean $b_i^\top \mu$ and variance $b_i^\top \Sigma b_i$. Moreover, the inequality system $B^{(i)}L^{(i)}y\le x^{(i)}-B^{(i)}w^{(i)}$ is nondegenerate.
\end{theorem}
Theorem \ref{thm:gradient} builds that assessing the gradient of the probability function can be reduced to computing the probability values. Based on the result, Henrion and M\"{o}ller (2012) developed some nonlinear optimization procedure to solve the CCP. Wei et al. (2024a) studied an appointment scheduling optimization problem. They required that the probability that the overtime of the operating room exceeds a threshold should be controlled within a confidence level, and formulated the problem to a CCP. They transformed the nonlinear cost function inequality within the probability function equivalently to a linear inequality system and reformulated the CCP as a JCCP where the constraint takes the form of (\ref{equ:lsep}). They considered a normal distribution and used the gradient expression of Henrion and M\"{o}ller (2012), and implemented a gradient based optimization algorithm to solve the CCP. Note that Henrion and M\"{o}ller (2012) focused on the linear cost functions within the chance constraint. The nonlinear setting was also studied in the literature. For instance, van Ackooij and Henrion (2014) derived the gradient formulas for some probability functions with nonlinear cost functions and Gaussian and Gaussian-like distributions. See also the studies reviewed therein.  

In some applications, it may be inappropriate to model the random vector as a normal distribution. Hu et al. (2022) studied using the Gaussian mixture model (GMM) to model the random distribution in the CCP. The density of a GMM takes the form
\begin{equation}\label{equ:mixGauss}
p_*(z)=\sum_{j=1}^K \pi_j \mathcal
{N}(\left.z\right|\mu_j,\Sigma_j),
\end{equation}
where $\mathcal
{N}(\left.z\right|\mu_j,\Sigma_j), j=1,\cdots,K$ are normal distributions and $\pi_j, j=1,\cdots,K$ are nonnegative weights summing to 1. The GMM possesses sufficient flexibility and can approximate a distribution to any accuracy under certain regularity conditions. It provides a unified input modeling approach for CCP. The GMM can be viewed as a probabilistic generative model. That indicates we can fit a GMM based on the information available and then use the distribution to generate more sample. The probabilistic generative model has been very popular in recent years. 

Suppose that the underlying distribution is a mixture distribution, with the expression
\begin{equation}\label{equ:mixp}
p_*(z)=\sum_{j=1}^K \pi_j p_j(z),
\end{equation}
where $p_j(z), j=1,\cdots,K$ are component distributions. Hu et al. (2022) showed that under (\ref{equ:mixp}),  
\begin{equation}\label{equ:GMMsep}
{\Pr}_{\sim P_*}\left\{c(x,\xi)\le 0\right\}=\sum_{j=1}^K
\pi_j{\Pr}_{\sim P_j}\left\{c(x,\xi)\le 0\right\},
\end{equation}
and 
\begin{equation}\label{equ:GMMgradsep}
\nabla_x {\Pr}_{\sim P_*}\left\{c(x,\xi)\le 0\right\}=\sum_{j=1}^K
\pi_j\nabla_x {\Pr}_{\sim P_j}\left\{c(x,\xi)\le 0\right\}.
\end{equation}
The equations build a bridge between the mixture distribution and the component distributions in a CCP, and is particularly useful for a GMM setting. The probability function under the normal distribution may have nice properties. Those properties may be generalized for the GMM. Given the flexibility and powerfulness of the GMM in modeling randomness, such generalization is potentially very useful.  

Hu et al. (2022) discussed several classes of CCP under the GMM (\ref{equ:mixGauss}). They considered the CCP with (\ref{equ:lsep}). Because of (\ref{equ:GMMsep})-(\ref{equ:GMMgradsep}), they showed that under GMM, the gradient of the probability function in (\ref{equ:lsep}) can be derived by using the gradient formula of the probability function under normal distribution provided in Henrion and M\"{o}ller (2012). Hu et al. (2022) also considered the following CCP
\begin{eqnarray}\label{eq:ccp:lp}
&\displaystyle\mathop{\mbox{minimize}}_{x\in X} & h(x)\\
&\mbox{ subject to} & {\Pr}_{\sim P_*}\left\{f_{0}(x)+f(x)^\top \xi
\le 0\right\}\ge 1-\alpha, \label{eq:ccp:lpc}
\end{eqnarray}
where $f(x)=\left(f_1(x),\cdots, f_k(x)\right)^\top$ and $f_i(x), i=0,1,\cdots,k$ are real valued functions
defined on $X$. Under the GMM, they reformulated the constraint (\ref{eq:ccp:lpc}) as 
 \begin{equation}\label{equ:normproduct}
\sum_{j=1}^K \pi_j \Phi\left(\frac{-\left[f_0(x)+f(x)^\top
\mu_j\right]}{\sqrt{f(x)^\top \Sigma_j f(x)}}\right) \ge 1-\alpha
\end{equation}
where $\Phi(\cdot)$ is the standard normal distribution function. Given that $\Phi(\cdot)$ can be evaluated readily, the constraint (\ref{equ:normproduct}) can be treated as a deterministic constraint.  They then suggested to use some nonlinear optimization procedure to solve the CCP with (\ref{equ:normproduct}). 

To achieve the global optimal value, Hu et al. (2022) proposed to use a branch-and-bound procedure to solve Problem (\ref{eq:ccp:lp}). In this approach, Hu et al. (2022) transformed the chance constraint (\ref{eq:ccp:lpc}) to the following set of constraints by bringing in a new decision vector $y=\left(y_1,\cdots,y_K\right)^\top$
\begin{equation}\label{equ:normproductsum}
{\Pr}_{\sim P_j}\left\{f_{0}(x)+f(x)^\top \xi \le 0\right\} \ge
y_j,  j=1,\cdots,K, \nonumber
\end{equation}
\begin{equation*}
\sum_{j=1}^K \pi_j y_j \ge 1-\alpha,\ 0\le y_j\le 1, j=1,\cdots, K.
\end{equation*}
They then reformulated Problem (\ref{eq:ccp:lp}) into the following optimization problem 
\begin{eqnarray}\label{eq:ccp:mipref}
&\displaystyle\mathop{\mbox{minimize}}_{x\in X, y} & h(x)\\
&\mbox{ subject to} & \Phi^{-1}(y_j)\sqrt{f(x)^\top \Sigma_j
f(x)}+\left(f_0(x)+f(x)^\top \mu_j\right)\le 0, j=1,\cdots, K, \label{eq:ccp:lpref1} \nonumber\\
&& \sum_{j=1}^K\pi_jy_j \ge 1-\alpha, \ 0\le y_j\le 1, j=1,\cdots, K. \nonumber \label{eq:ccp:lpref12}
\end{eqnarray}   
Problem (\ref{eq:ccp:mipref}) becomes a deterministic optimization problem given that $\Phi^{-1}$ can be evaluated readily. However, the decision variables $x$ and $y$ are coupled together, making the problem still difficult to address. Hu et al. (2022) studied a linear structure for the functions $f_{0}(x)$ and $f(x)$  and developed a spatial branch-and-bound algorithm to solve Problem (\ref{eq:ccp:mipref}) to the global optimality. Their idea is to partition the space of $y$ into sub-rectangles and relax Problem (\ref{eq:ccp:mipref})  within each sub-rectangle. Consider a sub-rectangle $\Delta=\{y\in \mathbb{R}^{K}\mid \underline{y}_j\le
y_j\le\overline{ y}_j,j=1,\ldots,K\}$ and the following problem 
\begin{eqnarray}\label{eq:ccp:miprefdelta}
&\displaystyle\mathop{\mbox{minimize}}_{x\in X, y\in\Delta} & h(x)\\
&\mbox{ subject to} & \Phi^{-1}(y_j)\sqrt{f(x)^\top \Sigma_j
f(x)}+\left(f_0(x)+f(x)^\top \mu_j\right)\le 0, j=1,\cdots, K, \label{eq:ccp:lprefdelta1} \nonumber \\
&& \sum_{j=1}^K\pi_jy_j \ge 1-\alpha. \nonumber \label{eq:ccp:lprefdelta2}
\end{eqnarray}
Hu et al. (2022) constructed the following problem 
\begin{eqnarray}\label{eq:ccp:miprefc1}
&\displaystyle\mathop{\mbox{minimize}}_{x\in X, y\in\Delta} & h(x)\\
&\mbox{ subject to} & \Phi^{-1}(\underline{y}_j)\sqrt{f(x)^\top \Sigma_j
f(x)}+\left(f_0(x)+f(x)^\top \mu_j\right)\le 0, j=1,\cdots, K, \label{eq:ccp:miprefc1c1} \nonumber\\
&& \sum_{j=1}^K\pi_jy_j \ge 1-\alpha.\nonumber \label{eq:ccp:miprefc1c2}
\end{eqnarray}
and showed that the optimal value of Problem (\ref{eq:ccp:miprefc1}) is a lower bound of the optimal value of Problem (\ref{eq:ccp:miprefdelta}). The lower bound is used to develop the spatial branch-and-bound algorithm. 

The basic idea of the  spatial branch-and-bound algorithm in Hu et al. (2022) is as follows. Let $v_{\mbox{opt}}$ denote the optimal value of Problem (\ref{eq:ccp:lp}). If solving Problem (\ref{eq:ccp:miprefc1}) generates an optimal solution $(x^*,y^*)$ such that $x^*$ is also feasible to Problem (\ref{eq:ccp:lp}), the optimal value of Problem (\ref{eq:ccp:miprefc1}) provides an upper bound of $v_{\mbox{opt}}$. Then the upper bound of $v_{\mbox{opt}}$ can be updated and the sub-rectangle $\Delta$ can be removed from the set, denoted as $\Omega$, of the sub-rectangles. If Problem  (\ref{eq:ccp:miprefc1}) is infeasible or its optimal value is not better than the current attained optimal value of Problem (\ref{eq:ccp:lp}), the sub-rectangle $\Delta$ can also be removed. If solving Problem (\ref{eq:ccp:miprefc1}) does not generate a feasible solution $x^*$ for Problem (\ref{eq:ccp:lp}), then the optimal value of Problem (\ref{eq:ccp:miprefc1}) provides a lower bound of the optimal value of Problem (\ref{eq:ccp:miprefdelta}). In this case, $\Delta$ remains in $\Omega$. After the computation, some of the inferior sub-rectangles in $\Omega$ are removed and the set $\Omega$ is updated. In the next iteration of the procedure, one sub-rectangle from $\Omega$ is selected based on some selection rule and is partitioned into smaller pieces based on some partition rule for further exploration. The process is repeated until the termination rule is met. Hu et al. (2022) provided the detailed algorithm and built properties of the algorithm. 

Wei et al. (2024b) studied how to speed up the spatial branch-and-bound algorithm proposed in Hu et al. (2022). They designed an enhanced pruning strategy and an augmented
branching strategy, and proposed an enhanced branch-and-bound algorithm. They conducted extensive numerical experiments to test the efficiency of the enhanced algorithm. Their numerical experiments also suggested that the nonlinear optimization algorithm could often obtain the global optimal value of Problem (\ref{eq:ccp:lp}) under the linear setting. Pang et al. (2024) studied a CCP with the form (\ref{eq:ccp}) where the underlying distribution is a GMM and the cost function $c(x,\xi)$ is a quadratic form of the random vector $\xi$
\begin{eqnarray*}
c(x, \xi)=\frac{1}{2}\xi^{\top} {\bf A}(x)\xi+{\bf a}(x)^{\top}\xi+a(x), 
\end{eqnarray*}
where ${\bf A}:\mathbb{R}^{n}\rightarrow\mathbb{R}^{m\times m}$, ${\bf a}:\mathbb{R}^{n}\rightarrow\mathbb{R}^{m}$ and $a:\mathbb{R}^{n}\rightarrow\mathbb{R}$ are mappings of appropriate dimensions. They showed that under certain conditions, the asymptotic distribution of $c(x, \xi)$ under the GMM is a univariate GMM. As a very interesting finding, they showed that under some conditions that are easy to satisfy, the chance constraint (\ref{eq:ccp:lpc}) with linear cost function actually defines a convex set. That explains why the nonlinear optimization procedure could often obtain the global optimum for Problem (\ref{eq:ccp:lp}) under the linear setting.

\section{Conclusions and Further Discussions}\label{sec:discussion} 

In this paper, we have discussed some of the studies on the computing of risk measures.  Specifically, we have discussed the computing of OCE risk measure, utility-based shortfall risk, and chance constrained program. A general remark is that the stochastic optimization theories and methods are important tools for addressing the problems of quantitative risk management. An interesting observation is that when using SAA to conduct estimation or optimization of risk measures, we can often achieve a convergence rate of order $n^{-1/2}$, which is the conventional convergence rate of the Monte Carlo simulation. Our studies only focused on some typical problems. Moreover, we note that this work has more focused on the author's studies and has not aimed to provide a comprehensive review of the research and advancement in this area. Surrounding the quantitative analysis of risk measures, there are a number of important topics and interesting directions we did not cover and explore. 

One important topic is the nested simulation of risk measures. In risk management, one often needs to assess the quantity 
\begin{equation*}
\alpha=\rho \left(\E\left[Y| X\right]\right), 
\end{equation*} 
where $X$ is a random vector, $\E\left[Y| X\right]$ is a conditional expectation of the random variable $Y$ given $X$, and $\alpha$ is the risk measure $\rho$ of the conditional expectation. The estimation problem is common in financial risk management and there have been deep studies on this topic, see, e.g., Gordy and Juneja (2010), Broadie et al. (2015), Hong et al. (2017), Zhang et al. (2022) and Wang et al. (2024). Liu and Zhang (2024) provided a comprehensive tutorial on nested simulation. Nested simulation for various risk measures is an important topic for study.   

The paper reviewed the statistical properties of various estimators for assessing or optimizing risk. However, many of the results are built in an asymptotic regime.  In practice, the data size or computational budget is often limited. It is beneficial to build some nonasymptotic results for the computation. As mentioned, Hegde et al. (2024) conducted such analysis for the SR. How to establish nonasymptotic statistical properties for the risk measures (e.g., derive some safe confidence intervals  for a quantity with a finite number of sample) could be an interesting research direction. 

Another important problem is regarding the risk and uncertainty. Note that the risk measure is a functional of the distribution of a random position. When one uses a risk measure, the uncertainty of the underlying distribution is also a significant concern faced by the decision maker. The input uncertainty is an important problem for decision under uncertainty. Hu and Hong (2022) proposed a robust simulation approach to address the input uncertainty in simulation. They studied robust simulation of various performance/risk measures, including expectation, probability, OCE risk measure and SR. The basic idea is to build the connections between the risk measures and the expectation. In practice, how to integrate the risk and uncertainty to address different application scenarios and different decision attitudes is an important problem and deserves further study. Regarding the optimization under uncertainty, distributionally robust optimization (DRO) has been a very popular and important approach to address the problem, see, e.g., the recent review in Kuhn et al. (2025). There have been fruitful studies on the DRO of the aforementioned risk measures in the literature, e.g., Zhu and Fukushima (2009) and Guo and Xu (2019). Due to the limited scope of our study, we do not explore and discuss the topic of DRO.

\section*{References}

\begin{hangref}

\item Alexander, S., T. F. Coleman, Y. Li. 2006. Minimizing CVaR and VaR for a
portfolio of derivatives. {\it Journal of Banking and Finance}, {\bf
30} 583-605.

\item Artzner, P., F. Delbaen, J.-M. Eber, D. Heath. 1999.
Coherent measures of risk. {\it Mathematical Finance}, {\bf 9}(3)
203-228.

\item Ben-Tal, A., M. Teboulle. 1986. Expected utility, penalty functions, and duality in stochastic nonlinear programming. {\it Management Science}, {\bf32}(11) 1445-1466.

\item Ben-Tal, A., M. Teboulle. 2007. An old-new concept of convex risk
measures: The optimized certainty equivalent. {\it Mathematical
Finance}, {\bf 17}(3) 449-476.

\item Broadie, M., Y. Du, C. C. Moallemi. 2015. Risk estimation via regression. {\it Operations Research} {\bf63}(5) 1077-1097.

\item Calafiore, G., M. C. Campi. 2005. Uncertain convex programs:
Randomized solutions and confidence levels. {\it Mathematical
Programming}, {\bf 102} 25-46.

\item Calafiore, G., M. C. Campi. 2006. The scenario approach to
robust control design. {\it IEEE Transactions on Automatic Control},
{\bf 51} 742-753.

\item Charnes, A., W. W. Cooper, G. H. Symonds. 1958. Cost
horizons and certainty equivalents: An approach to stochastic
programming of heating oil. {\it Management Science}, {\bf 4}
235-263.

\item Cui, Y., J. Liu, and J. S. Pang. 2022. Nonconvex and nonsmooth approaches for affine
chance-constrained stochastic programs. {\it Set-Valued and Variational Analysis}, {\bf30}(3) 1149-1211.

\item De Farias, D. P., B. Van Roy. 2004. On constraint sampling in the linear
programming approach to approximate dynamic programming. {\it
Mathematics of Operations Research}, {\bf 29}(3) 462-478.

\item Dunkel, J., S. Weber. 2010. Stochastic root finding and efficient
estimation of convex risk measures. {\it Operations Research}, {\bf
58}(5) 1505-1521.

\item Durrett, R. 2019. {\it Probability: Theory and Examples}.
Cambridge University Press.

\item Feng, G., G. Liu. 2016. Conditional Monte Carlo: A change-of-variables approach. arxiv preprint arxiv:1603.06378. 

\item F\"{o}llmer, H., A. Schied. 2002. Convex measures of risk and trading
constraints. {\it Finance and Stochastics}, {\bf 6} 429-447.

\item Frittelli, M., E. R. Gianin. 2002. Putting order in risk measures. {\it Journal of Banking and Finance}, {\bf26}(7) 1473-1486.

\item Glynn, P. W., L. Fan, M. C. Fu, J. Q. Hu, Y. Peng. 2020. Central limit theorems for estimated functions at estimated points. {\it Operations Research}, {\bf 68}(5) 1557-1563.

\item Glasserman, P. 2004. {\it Monte Carlo Methods in Financial
Engineering}, Springer, New York.

\item Gordy, M. B., S. Juneja. 2010. Nested simulation in portfolio risk measurement. {\it Management Science}, {\bf 56}(9) 1658-1673.

\item Guo, S., H. Xu. 2019. Distributionally robust shortfall risk optimization model and its approximation. {\it Mathematical Programming}, {\bf174}(1) 473-498.

\item Hamm, A. M., T. Salfeld, S. Weber. 2013. Stochastic
root finding for optimized certainty equivalents. {\it Proceedings
of the 2013 Winter Simulation Conference}, 922-932.

\item Hegde, V., A. S. Menon, L. A. Prashanth, K. Jagannathan. 2024. Online estimation and optimization of utility-based shortfall risk. {\it Mathematics of Operations Research}. Articles in Advance, 1-32.  

\item Henrion, R., A. M\"{o}ller. 2012. A gradient formula for linear
chance constraints under Gaussian distribution. {\it Mathematics of
Operations Research}, {\bf 37}(3) 475-488.

\item Hong, L. J. 2009. Estimating quantile sensitivities. {\it Operations
Research}, {\bf 57} 118-130.

\item Hong, L. J., Z. Hu., G. Liu. 2014. Monte Carlo methods for
value-at-risk and conditional value-at-risk: A review. {\it ACM
Transactions on Modeling and Computer Simulation}, {\bf 24}(4)
Article 22.

\item Hong, L. J., G. Jiang. 2017. Gradient and hessian of joint probability function with applications on chance-constrained programs. {\it Journal of the Operations Research Society of China}, {\bf5} 431-455.

\item Hong, L. J., S. Juneja, G. Liu. 2017. Kernel smoothing for nested estimation with application to portfolio risk measurement. {\it Operations Research}, {\bf65}(3) 657-673.

\item Hong, L. J., G. Liu. 2009. Simulating sensitivities of
conditional value-at-risk. {\it Management Science}, {\bf 55}
281-293.

\item Hong, L. J., Y. Yang, L. Zhang. 2011. Sequential convex approximations
to joint chance constrained programs: A Monte Carlo approach. {\it
Operations Research}, {\bf 59} 617-630.

\item Hu, Z., L. J. Hong. 2022. Robust simulation with likelihood-ratio constrained input uncertainty. {\it INFORMS Journal
on Computing}, {\bf 34}(4) 2350-2367.

\item Hu, Z., L. J. Hong, L. Zhang. 2013. A smooth Monte Carlo approach to joint
chance constrained program. {\it IIE Transactions}, {\bf 45}(7)
716-735.

\item Hu, Z., W. Sun, S. Zhu. 2022. Chance constrained programs with Gaussian mixture models. {\it IISE Transactions}, {\bf 54}(12) 1117-1130.

\item Hu, Z., D. Zhang. 2018. Utility‐based shortfall risk: Efficient computations via Monte Carlo. {\it Naval Research Logistics}, {\bf 65}(5) 378-392.

\item Huber, P. J., E. M. Ronchetti. 2011. {\it Robust Statistics}. John Wiley \& Sons. 

\item Jorion, P. 2006. {\it Value at Risk}, 2nd ed. McGraw-Hill, New York.

\item Krokhmal, P. A. 2007. Higher moment coherent risk measures. {\it Quantitative Finance}, {\bf 7}(4) 373-387. 

\item Kuhn, D., S. Shafiee, W. Wiesemann. 2025. Distributionally robust optimization. {\it Acta Numerica}, {\bf 34} 579-804. 

\item Küçükyavuz, S., R. Jiang. 2022. Chance-constrained optimization under limited distributional information: A review of reformulations based on sampling and distributional robustness. {\it EURO Journal on Computational Optimization}, {\bf10} 100030.

\item Liu, G., K. Zhang. 2024. A tutorial on nested simulation. {\it Proceedings of the 2024 Winter Simulation Conference}, 1-15. 

\item Luedtke, J., S. Ahmed. 2008. A sample approximation approach for
optimization with probabilistic constraints. {\it SIAM Journal on
Optimization}, {\bf 19}(2) 674-699.

\item Pagnoncelli, B. K., S. Ahmed, A. Shapiro. 2009. Sample average
approximation method for chance constrained programming: Theory and
applications. {\it Journal of Optimization Theory and Applications},
{\bf 142} 399-416.

\item Pang, X., Z. Hu, S. Zhu. 2024. Chance constrained program with quadratic randomness: A unified approach based on Gaussian mixture distribution. Working paper. 

\item Pe\~na-Ordieres, A., J. Luedtke, A. W\"{a}chter. 2020. Solving chance-constrained
problems via a smooth sample-based nonlinear approximation. {\it SIAM Journal on Optimization}, {\bf 30}(3) 2221-2250.

\item Rockafellar, R. T., J. O. Royset. 2015. Measures of residual risk with
connections to regression, risk tracking, surrogate models, and
ambiguity. {\it SIAM Journal on Optimization}, {\bf 25}(2)
1179-1208.

\item Rockafellar, R. T., S. Uryasev. 2000. Optimization of conditional
value-at-risk. {\it The Journal of Risk}, {\bf 2} 21-41.

\item Rockafellar, R. T., S. Uryasev. 2013. The fundamental risk quadrangle in
risk management, optimization and statistical estimation. {\it
Surveys in Operations Research and Management Science}, {\bf 18}
33-53.

\item Ruszczy\'nski, A., A. Shapiro. 2006. Optimization of convex risk
functions. {\it Mathematics of Operations Research}, {\bf 31}
433-452.

\item Shan, F., L. Zhang, X. Xiao. 2014. A smoothing function approach to joint chance-constrained programs. {\it Journal of Optimization Theory and Applications}, {\bf163} 181-199.

\item Shapiro, A., D. Dentcheva, A. Ruszczy\'nski. 2014. {\it Lectures on
Stochastic Programming: Modeling and Theory}, Second Edition. SIAM,
Philadelphia.

\item Trindade, A. A., S. Uryasev, A. Shapiro, G. Zrazhevsky. 2007.
Financial prediction with constrained tail risk. {\it Journal of
Banking and Finance}, {\bf 31} 3524-3538.

\item van Ackooij, W., R. Henrion. 2014. Gradient formulae for nonlinear probabilistic constraints with Gaussian and Gaussian-like distributions. {\it SIAM Journal on Optimization}, {\bf 24}(4) 1864-1889.

\item Wang, W., Y. Wang, X. Zhang. 2024. Smooth nested simulation: Bridging cubic and square root convergence rates in high dimensions. {\it Management Science}, {\bf 70}(2) 9031-9057.

\item Wei, J., Z. Hu, J. Luo. 2024a. Appointment scheduling optimization with chance constraints in a single-server consultation system. {\it Systems Engineering — Theory $\&$ Practice}, {\bf44}(10) 3400-3417. (In Chinese). 

\item Wei, J., Z. Hu, J. Luo, S. Zhu. 2024b. Enhanced branch-and-bound algorithm for chance constrained programs with Gaussian mixture models. {\it Annals of Operations Research}, {\bf 338}(2) 1283-1315.

\item Zhang, K., G. Liu, S. Wang. 2022. Bootstrap-based budget allocation for nested simulation. {\it Operations Research}, {\bf70}(2) 1128-1142.

\item Zhu, S., M. Fukushima. 2009. Worst-case conditional value-at-risk with application to robust portfolio management. {\it Operations Research}, {\bf 57}(5) 1155-1168.


\end{hangref}

\end{document}